\documentclass[12pt,a4paper]{amsart}
\usepackage{amssymb}
\usepackage{xypic}

\newtheorem{thm}{Theorem}[section]
\newtheorem{corol}[thm]{Corollary}

\newtheorem{prop}[thm]{Proposition}
\theoremstyle{definition}
\newtheorem{defin}[thm]{Definition}

\theoremstyle{remark}
\newtheorem{remark}[thm]{Remark}
\newtheorem{example}[thm]{Example}

\numberwithin{equation}{section}

\newcommand\Img{\operatorname{\mathcal Im}}

\newcommand\iso{\kern.35em{\raise3pt\hbox{$\sim$}\kern-1.1em\to}\kern.3em}

\newcommand\F{\mathbb F}
\newcommand\A{\mathbb A}
\newcommand\Pp{\mathbb P}
\newcommand\m{\mathfrak m}
\newcommand\Cc{\mathcal C}
\newcommand\Fc{\mathcal F}
\newcommand\Lcc{\mathcal L}
\newcommand\Oc{\mathcal O}
\newcommand\Div{\operatorname {Div}}
\newcommand\Res{\operatorname {Res}}
\newcommand\Spec{\operatorname {Spec}}
\newcommand\Proj{\operatorname {Proj}}

\begin{document}
\title{Convolutional Codes of Goppa Type}
\author{J.A Dom\'{\i}nguez P\'erez}
\author{J.M Mu\~noz Porras}
\author{G. Serrano Sotelo}
\email{jadoming@usal.es, jmp@usal.es and laina@usal.es}
\address{Departamento de Matem\'aticas, Universidad de Salamanca,
Plaza de la Merced 1-4, 37008 Salamanca, Spain}
%\date{\today}
\thanks {This research was partially supported by the Spanish DGESYC
through
research project BMF2000-1327 and by the ``Junta de Castilla y Le\'on''
through
research projects SA009/01 and SA032/02.}

\begin{abstract} A new kind of Convolutional Codes generalizing Goppa
Codes is proposed.
This provides a systematic method for constructing convolutional codes
with prefixed properties.
In particular, examples of Maximum-Distance Separable
(MDS) convolutional codes are obtained.
\end{abstract}

\maketitle

\section{Introduction}

The aim of this paper is to propose a definition of Convolutional
Goppa Codes (CGC). This definition will provide an
algebraic method for constructing Convolutional Codes with prescribed invariants.

We propose a definition of CGC in terms of families of curves $X\to
\A^1$ parametrized by the affine line
$\A^1= \Spec \F_q[z]$ over a finite field $\F_q$. In this setting,
the usual definition of a Goppa Code as the code
obtained by evaluation of sections at several rational points, is
translated as a code obtained by evaluation (of
sections of some invertible sheaf over $X$) along several sections of
the fibration $X\to \A^1$.

The paper is organized as follows.

In \S 2 we offer a summary on Goppa Codes following \cite{HovLP},
\cite{vLvG}, and using the standard notations of
Algebraic Geometry \cite{Har}.

\S 3 is devoted to giving the general definition of CGC and gives some
general results.

In \S 4 we study the case of a trivial fibration of projective lines
over $\A^1$ and we conclude giving some
explicit examples of MDS convolutional codes.

We freely use the standard notations of  abstract Algebraic
Geometry as can be found in \cite{Har}. After the
works of V. Lomadze \cite{Lom}, J. Rosenthal and R. Smarandache
\cite{RoSm}, \cite{SmRo}, there is evidence that
the use of methods of Algebraic Geometry can be relevant to the study
of Convolutional Codes. This paper is a step
in favor of that evidence.

\section{Background on Algebraic Geometry and Goppa Codes}

In this Section we summarize the basic definitions about Goppa Codes,
constructed using methods of Algebraic Geometry
  (see \cite{HovLP}, \cite{vLvG}).

Let $X$ be a geometrically irreducible, smooth and projective curve
over the finite field $\F_q$. Let $p_1,\dots,
p_n$ be $n$ different $\F_q$-rational points of $X$, and $D$ the
divisor $D=p_1+\dots+p_n$. Let $G$ be another
effective divisor with support disjoint from $D$. The Goppa code
$C(G,D)$ defined by $(G,D)$ is the linear code of
length $n$ over $\F_q$ defined as the image of the linear map
$$\aligned \alpha\colon  L(G)&\to \F_q^n\cr f&\mapsto
(f(p_1),\dots,f(p_n))\,,\endaligned$$
where $L(G)$ is the complete linear series defined by $G$. That is,
let $\F_q(X)$ be the field of rational functions
over the curve $X$,
$$L(G)=\left\{ f\in \F_q(X) \text{ such that } \Div(f)+G\geq 0\right\}\,.$$
The Goppa code has dimension
$$k=\dim C(G,D)=\dim L(G)-\dim L(G-D)\,.$$
Let $g$ be the genus of $X$; if we assume the inequality
$2g-2<\deg(G)<n$, then one has
$$k=\deg(G)-g+1\,,$$
and the minimum distance $d$ of $C(G,D)$ satisfies the inequality
$$d\geq n-\deg (G)\,.$$
Let $\Oc_X(D)$ be the invertible sheaf on $X$ defined by the divisor
$D$. One has the following exact sequence of
sheaves
$$0\to \Oc_X(-D)\to \Oc_X\to \Oc_D\to 0\,,$$
where $\Oc_D\simeq \Oc_{p_1}/\m_{p_1}\times \dots \times
\Oc_{p_n}/\m_{p_n}\simeq \F_q\times\overset {n)}\dots
\times\F_q$. Tensoring the above exact sequence by $\Oc_X(G)$, one obtains
$$0\to \Oc_X(G-D)\to \Oc_X(G)\to \Oc_D\to 0\,.$$
By taking global sections, we obtain an exact sequence of cohomology
$$\aligned
0\to H^0(X,\Oc_X(G-D))\to H^0(X,\Oc_X(G)) \overset\alpha\to \Oc_D\to
H^1&(X,\Oc_X(G-D))\to\cr
  &\to H^1(X,\Oc_X(G))\to 0\,,\endaligned$$
where $L(G)=H^0(X,\Oc_X(G))$ and $\alpha$ is the evaluation map defined above.

In the case $2g-2<\deg(G)<n$, one has the exact sequence
\begin{equation}\label{sq}
0\to H^0(X,\Oc_X(G)) \overset\alpha\to \Oc_D\to H^1(X,\Oc_X(G-D))\to 0\,.
\end{equation}

Let $\omega_X$ be the dualizing sheaf of $X$, which is isomorphic to
the sheaf of regular $1$-forms over $X$;
$H^0(X,\omega_X)$ is the $\F_q$-vector space of global regular
$1$-forms over $X$, which is of dimension $g=$ genus
of $X$.

By Serre's duality (\cite{Har}), there exist canonical isomorphisms of
$\F_q$-vector spaces
$$H^1(X,\Lcc)^\ast\simeq H^0(X,\omega_X\otimes \Lcc^{-1})$$
for every invertible sheaf $\Lcc$ on $X$. Given a divisor $D$ over
$X$, we shall denote by $\Omega(D)$ the vector
space $H^0(X,\omega_X\otimes \Oc_X(-D))$.

The dual Goppa code, $C^\ast(G,D)$, associated with the Goppa code
$C(G,D)$ is defined as the linear code of length
$n$ over $\F_q$ given by the image of the linear map
$$\aligned \alpha^\ast\colon  \Omega(G-D)&\to \F_q^n\cr \eta &\mapsto
(\Res_{p_1}(\eta),\dots,\Res_{p_n}(\eta))\,,\endaligned$$
Let us take duals in the exact sequence (\ref{sq}):
$$
0\to H^1(X,\Oc_X(G-D))^\ast \overset\beta\to \Oc_D^\ast
\overset{\alpha^t}\to H^0(X,\Oc_X(G))^\ast\to 0\,.
$$
By Serre's duality, one has isomorphisms
$$\aligned
&H^1(X,\Oc_X(G-D))^\ast\simeq \Omega(G-D)\,,\cr
&H^0(X,\Oc_X(G))^\ast\simeq H^1(X,\omega_X\otimes \Oc_X(-G))\,,
\endaligned$$
and the above sequence is the  cohomology sequence
induced by the exact sequence of sheaves
$$
0\to \omega_X(-G) \to \omega_X(D-G)\to
\omega_X(D-G)\otimes_{\Oc_X}\Oc_D\to 0\,,
$$
where we denote $\omega_X(-G)=\omega_X\otimes \Oc_X(-G)$, and
$\beta$ is precisely the map $\alpha^\ast$
defining $C^\ast(G,D)$.

Given a linear series $\Gamma\subseteq H^0(X,\Oc_X(G))$,
that is, a  vector subspace defining a family of divisors
linearly equivalent to $G$, we define the Goppa code $C(\Gamma,D)$
associated whit $\Gamma$  and $D$ as the image of
the homomorphism $\alpha_{\vert_\Gamma}$:
$$\xymatrix@R2pt{H^0(X,\Oc_X(G)) \ar[r]^-{\alpha} & \Oc_D \\
\bigcup\vert & \\
\Gamma \ar[uur]_>>>>>>>>>>{\alpha_{\vert_\Gamma}}}$$
When $\Gamma\varsubsetneqq H^0(X,\Oc_X(G))$, we shall say that
$C(\Gamma,D)$ is a non-complete Goppa code.

\section{Convolutional Goppa Codes}

We shall
contruct a kind of convolutional code  that generalizes the notion of
Goppa codes. These codes will be associated with
families of algebraic curves.

Given an algebraic variety $S$ over $\F_q$, a family of projective
algebraic curves
parametrized by $S$ is a morphism of algebraic varieties
$\pi\colon X\to S$,
such that $\pi$ is a projective and flat morphism whose fibres
$X_s=\pi^{-1}(s)$ are smooth and geometrically
irreducible curves over $\F_q(s)$ (the residue field of $s\in S$).

Let us consider a family of curves $X\overset \pi
\to U$ parametrized by $U=\Spec \F_q[z]=\A^1$. Given a closed point
$u\in U$ with residue field $\F_q(u)$, the fibre
$X_u=\pi^{-1}(u)$ is a curve over the finite field  $\F_q(u)$.

Let $p_i$, $1\leq i\leq n$, be $n$ different sections, $p_i\colon U\to
X$, of the projection $\pi$. These sections
define a Cartier divisor on $X$:
$$D=p_1(U)+\dots+p_n(U)\,,$$
which is flat of degree $n$ over the base $U$ (\cite{Har}).

Note that given a coherent sheaf $\Fc$ on $X$, the
cohomology groups $H^i(X,\Fc)$ are finite $\F_q[z]$-modules
and $H^i(X,\Fc)=0$ for $i\geq 0$ (see \cite{Har} III).

Let $\Lcc$ be an invertible sheaf over $X$. One has an exact sequence
of sheaves on $X$
\begin{equation}\label{ls}
0\to \Lcc(-D) \to \Lcc\to \Oc_D\to 0\,,
\end{equation}
which induces a long exact cohomology sequence
\begin{equation}\label{pls}
0\to H^0(X,\Lcc(-D)) \to H^0(X,\Lcc)\overset\alpha\to H^0(X,\Oc_D)\to
H^1(X,\Lcc(-D))\to
H^1(X,\Lcc)\to 0\,.
\end{equation}

Let $r$ be the degree of $\Lcc$ in each fibre of $\pi$ (which is
independent of the fibre) and let $g$ be the genus of
any fibre of $\pi$ (also independent of the fibres).

\begin{prop}\label{cond} Let us assume that $2g-2<r$. Then, one has
that $H^1(X,\Lcc)=0$ and $H^0(X,\Lcc)$ is a free
$\F_q[z]$-module of rank $r-g+1$
\end{prop}
\begin{proof} Under the condition $2g-2<r$, one has that
$H^1(X_u,{\Lcc}_{\vert_{X_u}})=0$ for every
point $u\in U$. Note that
$H^i(X,\Fc)\widetilde{\ } =R^i\pi_\ast \Fc$ for every coherent
sheaf
$\Fc$ on
$X$ (\cite{Har} III), and applying (\cite{Har} III Corollary 12.9)
one concludes the proof.
\end{proof}

Under the hypothesis of Proposition \ref{cond}, there exists an exact
sequence of $\F_q[z]$-modules
\begin{equation}\label{mpls}
0\to H^0(X,\Lcc(-D)) \to H^0(X,\Lcc)\overset\alpha\to H^0(X,\Oc_D)\to
H^1(X,\Lcc(-D))\to
0\,.
\end{equation}
where $H^0(X,\Oc_D)$ is a free $\F_q[z]$-module of rank $n$.

\begin{remark}
Let $\eta\in U$ be the generic point of $U$, whose residue field is
$\F_q(z)$; the fibre  $X_\eta=\pi^{-1}(\eta)$ is
a smooth, irreducible curve over $\F_q(z)$.
Note that $p_1(\eta),\dots, p_n(\eta)$ are $n$ different
$\F_q(z)$-rational points of the curve $X_\eta$.
One then has a canonical decomposition of $H^0(X,\Oc_D)_\eta$ as a
$\F_q(z)$-algebra
$$H^0(X,\Oc_D)_\eta=\F_q(z)\times\overset{n)}\dots \times \F_q(z)\,.$$
\end{remark}

Given a $\F_q[z]$-module $M$, let us denote by $M_\eta$ the
$\F_q(z)$-vector space
$$M_\eta=M\otimes_{{\F}_q[z]}\F_q(z)\,.$$
The sequence (\ref{mpls}) induces an exact sequence of $\F_q(z)$-vector spaces
\begin{equation}\label{mmpls}
0\to H^0(X,\Lcc(-D))_\eta \to
H^0(X,\Lcc)_\eta\overset{\alpha_\eta}\to H^0(X,\Oc_D)_\eta\to
H^1(X,\Lcc(-D))_\eta\to
0\,.\end{equation}

\begin{defin}\label{cgc} The complete convolutional  Goppa code
associated with $\Lcc$ and $D$ is the image of the
homomorphism $\alpha_\eta$
$$\Cc(\Lcc,D)=\Img\left(
H^0(X,\Lcc)_\eta\overset{\alpha_\eta}\longrightarrow
H^0(X,\Oc_D)_\eta \simeq
\F_q(z)^n\right)\,.$$
Given a free submodule $\Gamma\subseteq H^0(X,\Lcc)$, the convolutional
Goppa code associated with $\Gamma$ and $D$ is
the image of ${\alpha_\eta}_{\vert_{\Gamma_\eta}}$
$$\Cc(\Gamma,D)=\Img\left( \Gamma_\eta \overset{\alpha_\eta}\longrightarrow
\F_q(z)^n\right)\,.$$
\end{defin}

\begin{remark}\label{rem} We use definition 2.4 of
\cite{McE} as definition of convolutional codes. Any
matrix defining $\alpha_\eta$ (respectively
${\alpha_\eta}_{\vert_{\Gamma_\eta}}$) is a generator matrix of
rational
functions for the code $\Cc(\Lcc,D)$ (resp. $\Cc(\Gamma,D)$).
\end{remark}

The canonical decomposition $H^0(X,\Oc_D)_\eta\simeq \F_q(z)^n$ as
$\F_q(z)$-algebras does not extend (in general) to
a decomposition $H^0(X,\Oc_D)\simeq \F_q[z]^n$ as rings.
In fact, one has a canonical isomorphism of rings
$H^0(X,\Oc_D)\overset\phi\iso\F_q[z]^n$
only when $p_1(U),\dots,p_n(U)$ are disjoint sections.
However,
$H^0(X,\Oc_D)$ is a free $\F_q[z]$-module; then, there exist
(non-canonical) isomorphisms of $\F_q[z]$-modules:
$$H^0(X,\Oc_D)\overset\phi\iso\F_q[z]\oplus\overset{n)}\dots \oplus
\F_q[z]\,,$$
which are not (in general) isomorphism of rings.

This allows us to give another definition of
convolutional Goppa codes.

\begin{defin}\label{cgct} Given a trivialization $\phi\colon
H^0(X,\Oc_D)\iso \F_q[z]^n$ as $\F_q[z]$-modules, one
defines the convolutional Goppa code $\Cc(\Lcc,D,\phi)$ as the image
of $\phi\circ\alpha$
$$H^0(X,\Lcc)\overset\alpha\to H^0(X,\Oc_D)\overset \phi\iso \F_q[z]^n\,.$$
Anagously, one defines the convolutional Goppa code $\Cc(\Gamma,D,\phi)$.
\end{defin}

Let us assume (for the rest of the paper) that the invariants
$(r,n,g)$ satisfy the inequality
$$2g-2<r<n\,.$$

\begin{prop}\label{suc} Under the above conditions on $(r,n,g)$,
$H^0(X,\Lcc(-D))=0$ and $H^1(X,\Lcc(-D))$ is a free
$\F_q[z]$-module. The following exact sequence is exact
\begin{equation}\label{corta}
0\to H^0(X,\Lcc)\overset{\alpha}\to H^0(X,\Oc_D)\to H^1(X,\Lcc(-D))\to 0\,.
\end{equation}
and remains exact when we take fibres over every point $u\in U$.
\end{prop}
\begin{proof} If $2g-2<r<n$,
$H^0(X_u,\Lcc(-D)_{\vert_{X_u}})=0$ for every point $u\in U$; and
applying
(\cite{Har} III Corollary 12.9) one concludes.
\end{proof}

\begin{corol} The convolutional code $\Cc(\Lcc,D,\phi)$ has dimension $k=r-g+1$ and length
$n$.
Every matrix defining $\phi\circ \alpha$ is a basic generator matrix
\cite{McE} for
$\Cc(\Lcc,D,\phi)$.
\end{corol}
\begin{proof} This is a direct consecuence of the last statement of
Proposition \ref{suc} and the characterization of
basic generator matrices of \cite{McE}.
\end{proof}

Let us consider the convolutional Goppa code $\Cc(\Gamma,D,\phi)$
defined by a submodule $\Gamma\subseteq H^0(X,\Lcc)$
and a trivilization $\phi$. With the above restrictions, one has:

\begin{prop}\label{prop} Every matrix defining
${\phi\circ\alpha}_{\vert_{\Gamma}}$ is a basic generator matrix for the code
$\Cc(\Gamma,D,\phi)$ if and only if $H^0(X,\Lcc)/\Gamma$ is a
torsion-free $\F_q[z]$-module.
\end{prop}
\begin{proof} The sequence (\ref{corta}) induces a diagram
$$\xymatrix@R12pt{  & 0 \ar[d] & 0 \ar[d] &   & \\
0 \ar[r] & \Gamma \ar[d]\ar[r]^-{\alpha_{\vert_{\Gamma}}} &
H^0(X,\Oc_D) \ar@{=}[d]\ar[r] &
H^1(X,\Gamma) \ar[d]\ar[r] &0\\
0 \ar[r] & H^0(X,\Lcc) \ar[d]\ar[r] & H^0(X,\Oc_D) \ar[r] &
H^1(X,\Lcc(-D)) \ar[d]\ar[r] &0\\
   & H^0(X,\Lcc)/\Gamma   &     &  0 & }$$
Then, the kernel of $H^1(X,\Gamma)\to H^1(X,\Lcc(-D))$ is isomorphic
to $H^0(X,\Lcc)/\Gamma$ and $H^1(X,\Lcc(-D))$ is
free. This implies that the torsion elements of $H^1(X,\Gamma)$ are
contained in $H^0(X,\Lcc)/\Gamma$, from which one
concludes the proof.
\end{proof}

The above results allow us to construct basic generator matrices for
the codes $\Cc(\Gamma, D,\phi)$. If
$p_1(U),\dots,p_n(U)$ are disjoint sections and $\phi$ the canonical
trivialization, this gives us a basic
generator matrix for
$\Cc(\Gamma, D)$. However, in general the codes $\Cc(\Gamma, D)$ and
$\Cc(\Gamma, D, \phi)$ are different.

Let us describe a geometric way to obtain a basic generator matrix
for $\Cc(\Lcc,D)$ and $\Cc(\Gamma, D)$.

Assume that the curves $p_1(U),\dots,p_n(U)$ meet transversally at
some points, and let $\bar X$ be the
blowing-up \cite{Har} of $X$ at these points. One has morphisms
$$\xymatrix{\bar X \ar[r]^-{\beta}\ar[dr]_-{\bar\pi=\pi\circ\beta} &
X\ar[d]^-\pi\\
  & U}$$
such that the proper transform of $D$ under $\pi$ is a divisor $\bar
D\subset \bar X$ satisfying
$$\bar D=p_1(U)\amalg\dots \amalg p_n(U)\overset \beta \to D\,,$$
and one has a canonical homomorphism of rings
$$0\to \Oc_D\to\beta_\ast \Oc_{\bar D}$$
which induces
$$0\to \pi_\ast\Oc_D\overset\beta\to\bar\pi_\ast \Oc_{\bar D}\simeq
\overset\sim{\F_q[z]^n}\,,$$
where $\bar\pi_\ast \Oc_{\bar D}\simeq \overset\sim{\F_q[z]^n}$ is
the canonical isomorphism of sheaves of rings.

$\beta^\ast\Lcc$ is an invertible sheaf on $\bar X$ and there exists a
canonical homomorphism
$$\beta^\ast\Lcc \to \Oc_{\bar D}\to 0\,,$$
whose kernel is $(\beta^\ast\Lcc)(-\bar D)$. This induces
$$0\to\Lcc\to\beta_\ast\beta^\ast\Lcc\to\beta_\ast\Oc_{\bar D}\,,$$
and taking global sections one obtains
$$0\to H^0(X,\Lcc)\overset\gamma\to
H^0(X,\beta_\ast\beta^\ast\Lcc)\overset\mu\to \F_q[z]^n\,.$$
The image of $\mu$ is precisely a free submodule of $\F_q[z]^n$ that
defines a basic generator matrix for
$\Cc(\Lcc,D)$.

Let us consider the sequence of homomorphisms
$$0\to H^0(X,\Lcc)\overset\alpha\to
H^0(X,\Oc_D)\overset\beta\hookrightarrow H^0(X,\Oc_{\bar
D})=\F_q[z]^n\,.$$
$\beta\circ\alpha$ is not in general a basic  matrix, since
$H^0(X,\Oc_{\bar D})/H^0(X,\Oc_D)$ has torsion. Let
us define
$$\bar H^0(X,\Lcc)=\{ p\in \F_q[z]^n \text{ such that } \lambda p\in
H^0(X,\Lcc)\text{ for some } \lambda\in \F_q[z]
\}\,.$$
$\bar H^0(X,\Lcc)/H^0(X,\Lcc)$ is a torsion module and
$\F_q[z]^n/\bar H^0(X,\Lcc)$ is torsion-free. Then, every matrix
defining the homomorphism $\bar H^0(X,\Lcc)\hookrightarrow \F_q[z]^n$
is a basic generator matrix for $\Cc(\Lcc,D)$.

This is an algebraic-geometric interpretation of Forney's
construction of the basic matrices of a convolutional
code \cite{For}.

\section{Convolutional Goppa Codes associated with the projective line}

Let $\Pp^1=\Proj \F_q[x_0,x_1]$ be the projective line over $\F_q$, and
$$X=\Pp^1\times U\overset \pi\to U=\Spec \F_q[z]$$
the trivial fibration. Let us denote by $t=x_1/x_0$ the affine
coordinate in $\Pp^1$, and by $p_\infty$ its infinity
point. Let us consider the following $n$ different sections of $\pi$
$$p_i\colon U\to \Pp^1\times U$$
defined in the coordinates $(t,z)$ by
$$p_i(z)=(\alpha_i z + \beta_i, z)\,,\ \alpha_i,\beta_i\in\F_q\,.$$
Let $D=p_1(U)+\dots+p_n(U)$ and let $\Lcc$ be the invertible sheaf on $X$
$$\Lcc=\pi_1^\ast \Oc_{\Pp^1}(r p_\infty)\otimes_{\F_q}\Oc_U\,,\ r<n\,,$$
The exact sequence (\ref{corta}) is in this case:
$$\xymatrix@R2pt{0 \to H^0(X,\Lcc) \ar[r]^-{\alpha} & H^0(X,\Oc_D)
\ar[r] & H^1(X,\Lcc(-D))\ar[r] & 0\,.\\
\ \ \ \ \ \ \vert\vert & \vert\vert\\
  H^0(\Pp^1,\Oc_{\Pp^1}(r p_\infty))\otimes\F_q[z] \ar[r]^-{\alpha} &
\F_q[z]^n}$$
Taking the fibres over the generic point $\eta$, and the canonical
trivialization
$\left(\pi_\ast\Oc_D\right)_\eta\simeq \F_q(z)^n$, the homomorphism
$\alpha_\eta$ is the evaluation map at the points
$p_1(\eta),\dots,p_n(\eta)$
$$\aligned \alpha_\eta&\colon H^0(\Pp^1,\Oc_{\Pp^1}(r
p_\infty))\otimes_{\F_q}\F_q(z)\to \F_q(z)^n\cr
\alpha_\eta&(t^j)=\left(t^j(p_1(\eta)),\dots,t^j(p_n(\eta))\right)=
\left((\alpha_1 z+\beta_1)^j,\dots,(\alpha_n z+\beta_n)^j\right)\,,
\endaligned$$
where $\{ 1,t,\dots, t^r\}$ is the ``canonical''  basis of
$H^0(\Pp^1,\Oc_{\Pp^1}(r p_\infty))$ in the affine
coordinate $t$. The convolutional code $\Cc(\Lcc,D)$ is a kind of
{\em generalized Reed-Solomon (RS)
code\/} (for $z=0$ we  obtain a classical RS-code).

Let $\Gamma\subseteq H^0(\Pp^1,\Oc_{\Pp^1}(r p_\infty))$ be the
linear subspace generated by $\{ t^s,\dots, t^r\}$.
The convolutional Goppa code $\Cc(\Gamma,D)$ is the image of the homomorphism
$$\aligned \alpha_\eta\colon\Gamma\otimes_{\F_q}&\F_q(z)\to \F_q(z)^n\cr
t^j&\longmapsto \alpha_\eta(t^j)\,,\ \text{ for\ } s\leq j\leq r\,.
\endaligned$$
In this case
$H^0(X,\Lcc)/\Gamma \simeq \left( H^0(\Pp^1,\Oc_{\Pp^1}(r p_\infty))
/\Gamma\right)\otimes_{\F_q} \F_q[z]$ is torsion-free.
Then, by Proposition \ref{prop} every matrix defining
$$\alpha\colon \Gamma\otimes_{\F_q}\F_q[z]\to H^0(X,\Oc_D)$$
is a basic generator matrix. To compute a matrix for
$\alpha$ explicitly, we need to fix an isomorphism of
$\F_q[z]$-modules
$$H^0(X,\Oc_D)\overset\phi\to\F_q[z]^n\,,$$
and this gives a generator matrix for $\Cc(\Gamma,D,\phi)$. However,
it would be desirable to compute basic matrices for
the codes $\Cc(\Gamma,D)$. We shall do this in general in a forthcoming paper.
Here we shall offer some explicit examples.

\begin{example}\label{ex} Let $a,b\in \F_q$ be two different non-zero
elements, and
$$p_i(z)=(a^{i-1} z +
b^{i-1}, z)\,, i=1,\dots, n\,, \text{ with $n<q$}\,. $$
The evaluation map ${\alpha_\eta}$ over $\Gamma$ is defined
by the matrix
\begin{equation}\label{matrix}
\left( \begin{matrix}
(z+1)^s  & (a z + b)^s & (a^2 z + b^2)^s &\dots & (a^{n-1}z+b^{n-1})^s\cr
(z+1)^{s+1}  & (a z + b)^{s+1} & (a^2 z + b^2)^{s+1} &\dots &
(a^{n-1}z+b^{n-1})^{s+1}\cr
\vdots  & \vdots & \vdots &\ddots & \vdots\cr
(z+1)^r  & (a z + b)^r & (a^2 z + b^2)^r &\dots & (a^{n-1}z+b^{n-1})^r\cr
\end{matrix}\right)\,.
\end{equation}
This matrix is a generator matrix for the code $\Cc(\Gamma,D)$. Using
this construction we can give concrete examples
of CGC of dimension $k=r-s+1$ that are  Maximum-Distance Separable
(MDS) convolutional codes, i.e., whose {\em free
distance\/}  attains the generalized Singleton bound \cite{RoSm}.

\bigskip

\begin{itemize}

\item If $s=r$, the convolutional Goppa code $\Cc(\Gamma,D)$
has  dimension $1$, degree $r$, and
(\ref{matrix}) is a
{\em canonical\/} (reduced and basic \cite{McE}) generator matrix. We
can list a few examples, where $k/n$, $\delta$ and $d$ are respectively the rate, the degree and the free
distance of the code.

\vbox to 0.5cm{\vfill}
\begin{center}
\begin{tabular}{|c|c|c|c|c|}
\hline
     &   &  &   &\\
    \em{field} &  \em{canonical generator matrix} &  $k/n$ &  $\delta$ &
$d$
\\
     &   &  &  & \\
\hline
     &   &  &   &\\
$\F_3=\{0,1,2\}$ &  $\left( \begin{matrix} z+1 & z+2 \end{matrix} \right)$  &  1/2  & 1 &   4  \\
     &   &  &   &\\
\hline
     &   &  &   &\\
$\F_4=\{0,1,\alpha,\alpha^2\}$ &  $\left( \begin{matrix} z+1 & z+\alpha & z+\alpha^2\end{matrix} \right)$  &  1/3  & 1 &  6  \\
   {\small where $\alpha^2+\alpha+1=0$ } &   &  &  & \\
\hline
     &   &  &  & \\
$\F_5=\{0,1,2,3,4\}$ &  $\left( \begin{matrix} (z+1)^2 & (z+2)^2 &
(z+4)^2\end{matrix}
\right)$  &  1/3  & 2 &   9  \\
     &   &  &   &\\
\hline
\end{tabular}
\end{center}
\bigskip

\noindent In these examples the sections
$p_1,\dots, p_n$ are disjoint, such that
one has $\Cc(\Gamma, D)=\Cc(\Gamma, D, \phi)$, where $\phi\colon
H^0(X,\Oc_D) \iso \F_q[z]^n$ is the
corresponding canonical trivialization.

\medskip

\item If $s<r$, let us take $a\in\F_q$  as a primitive
element.

Now, the matrix (\ref{matrix}) is
reduced, since the matrix of
highest-degree terms in each row is a  Vandermonde
matrix of rank $k$.
The sections $p_1,\dots, p_n$ are not disjoint, but
in some cases the matrix (\ref{matrix}) is actually basic
and we do not have to find an isomorphism of
$\F_q[z]$-modules,
$\phi\colon H^0(X,\Oc_D) \iso \F_q[z]^n$, in order to compute a
basic generator matrix for the code $\Cc(\Gamma,
D)$.

We present two examples of this situation.

\medskip
\begin{center}
\begin{tabular}{|c|c|c|c|c|}
\hline
     &   &  &  & \\
   \em{field}   &  \em{canonical generator matrix} &  $k/n$ & $\delta$ & $d$ \\
     &   &  &  & \\
\hline
     &   &  &   &\\
$\F_4$ &  $\left( \begin{matrix} 1 & 1 & 1\cr z+1 & \alpha
z+\alpha^2 & \alpha^2 z+\alpha\end{matrix}
\right)$  &  2/3  & 1 &  3  \\
     &   &  &  & \\
\hline
     &   &  &  & \\
$\F_5$ &  $\left( \begin{matrix} z+1 & 2z+3 & 4z+4& 3z+2\cr (z+1)^2 &
(2z+3)^2 & (4z+4)^2 & (3z+2)^2\end{matrix} \right)$  &
1/2  & 3 &  8
\\
     &   &  &  & \\
\hline
\end{tabular}
\end{center}
\medskip

\end{itemize}

\end{example}

\medskip
\noindent
{\bf Acknowledgments.} We thank F.J. Plaza Mart\'\i{}n and E. G\'omez
Gonz\'alez for many
enlightening comments, and J. Prada Blanco and J.I. Iglesias Curto
for helpful questions that helped us to
improve this paper.

\end{document}